\nonstopmode 
\magnification=\magstep 1 
\input amstex

\documentstyle{amsppt}
\loadbold \loadmsbm\loadeufm 
\topmatter
\title POINTS IN GENERIC POSITION AND CONDUCTOR OF VARIETIES WITH
ORDINARY MULTIPLE SUBVARIETIES OF CODIMENSION ONE\endtitle
\author Ferruccio Orecchia\endauthor
\address Ferruccio Orecchia, Dipartimento di Matematica e Applicazioni "R. Caccioppoli", Complesso 
Universitario di Monte S. Angelo
- Via Cintia , 80126 Napoli-Italy \endaddress
\email ORECCHIA\@MATNA2.DMA.UNINA.IT.\endemail\
\thanks Work partially supported by 
MURST\endthanks 

\abstract We extend  results of our previous papers, 
on ordinary multiple points of curves
[9], and on the computation of their conductor [8], 
to ordinary multiple subvarieties of codimension one. 
\endabstract
\endtopmatter

\document
\head 1. Introduction.\endhead 
 Let $A$ be the local ring, at a multiple (that is singular) point $x$,
 of an algebraic reduced curve $C$
with embedding dimension $emdim(A)=r+1$ and 
multiplicity $e(A)=e$. Let  $\frak m$ be the  maximal ideal of $A$.
$Spec(G(A))$ is the tangent cone and $Proj(G(A))$ the projectivized tangent cone
to $C$ at $x$. 

The scheme $Proj(G(A))$ is reduced   if and only if it 
consists of $e$ points of $\Bbb P^r$, that is,  the
tangent cone considered as a set consists of e lines of $\Bbb  A^{r+1}$ through $x$
(the tangents to the curve at $x$). In this case we say that   
$x$ is  an {\sl ordinary multiple point} (or an {\sl ordinary singularity}) of $C$
 [9, Lemma-Definition 2.1].  
Clearly if  
$Spec(G(A))$ is reduced then $Proj(G(A))$ is reduced.
The converse, in general, doesn't hold (see [9, Example 1] or [7, Section 4]), but 
if  $Proj(G(A))$  is reduced and consists of points  in generic $e$ position, 
then  $Spec(G(A))$ is reduced [9, Theorem 3.3].
Furthermore, if 
the points of $Proj(G(A))$ are in generic $e-1,e$
position, then  the conductor of $A$
in its normalization $\overline A$ is $\frak m^\nu$ where 
$\nu=Min\{n\in \Bbb N\mid e\leq $$n+r\choose r$$\}$ [8, Theorem 4.4].

In this paper first we extend these results to any one dimensional reduced local 
 ring
 $B$ with finite normalization $\overline B$. Let $\frak m$ be the
maximal ideal of $B$ and $K$
be the algebraic closure of the residue field $k(\frak m)$ of $B$. Set 
 $e=e(B)$, $emdim(B)=r+1$. Considering
 the ring $G(B)\otimes_{k(\frak m)} K$, instead of $G(A)$, 
we prove that, if $Proj(G(B) \otimes_{k(\frak m)}  K)$ is reduced and consists of 
points in
generic $e$ position, then $Spec(G(B) \otimes_{k(\frak m)} K)$ and $Spec(G(B))$
 are reduced [Theorem 5]. If 
in addition the points of $Proj(G(B) \otimes_{k(\frak m)}  K)$ are in generic $e-1$ position
then   the conductor of $B$
in $\overline B$ is $\frak m^\nu$ where 
$\nu=Min\{n\in \Bbb N\mid e\leq $$n+r\choose r$$\}$ [Theorem 11].\medskip      
Then we apply the  previous results to the ring $B=A_\frak p$,
where $\frak p$ is a prime ideal of codimension one of  a local reduced ring $A$, 
with finite
normalization $\overline A$.

Using also the properties of normal flatness
 we get the following result
on the conductor $\frak b$ of $A$ (in $\overline A$).\medskip
Let $e(A_\frak p)=e>1$,
 $emdim(A_\frak p)=r+1$ and $K$ be the algebraic closure of the residue field 
$k(\frak p)=A_\frak p/\frak pA_\frak p$ of $A$ at $\frak p$.
Assume $A/\frak p$  
regular, $\sqrt{ \frak b}=\frak p$ and $Proj(G(A_\frak p) \otimes_{k(\frak p)}  K)$  reduced
with   points in generic $e-1,e$ position. Given the the following  conditions:\medskip 
(a) $A$ is Cohen-Macaulay, 
 $emdim(A)=emdim(A_\frak p)+dim(A/\frak p)$ and $e(A)=e(A_\frak p)$;  \medskip
(b) $A$ is $S_2$ and  normally flat 
along $\frak p$;\medskip
(c) $\frak b$ is primary  and $A$ is normally flat 
along $\frak p$;\medskip
(d) $\frak b=\frak p^\nu$, where $\nu=Min\{n\in \Bbb N\mid e\leq $$n+r\choose r$$\}$;
\medskip 
then $(a) \Rightarrow (b) \Rightarrow (c) \Rightarrow (d)$ [Theorem 15].\medskip
In  particular  we get  the following statement.\medskip
Let $\frak p$ be a prime ideal of codimension one 
of a reduced
Cohen-Macaulay ring $A$ with finite normalization $\overline A$. Assume 
$A/\frak p$   
regular, $\sqrt{ \frak b}=\frak p$ and $Proj(G(A_\frak p) \otimes_{k(\frak p)}  K)$  
 reduced. 
If $emdim(A)=dim(A)+1$ and $e(A)=e(A_\frak p)=e$
then $\frak b=\frak p^{e-1}$ [Corollary 16].\medskip
The previous results have the following  geometrical consequences.\medskip
Let $X=Spec(R)$ be an algebraic variety and $Y=Spec(R/\frak q)$ be
 an irreducible codimension one subvariety   of $X$. Suppose that $e(R_\frak q)=e>1$
that is $Y$ is a {\sl multiple subvariety} of $X$, of multiplicity $e$. Let   
$emdim(R_\frak q)=r+1$ and $K$ be the algebraic closure of the residue field 
$k(\frak q)$ of $R$ in $\frak q$.
If $Proj(G(R_\frak q) \otimes_{k(\frak q)} K)$ is reduced and its points
are in generic $e-1,e$ position, then there exists an open nonempty subset $U$ of $Y$ 
such that
for every closed point $x$ of $U$  
the conductor of the local ring $A$, of $X$ at $x$,
 is  $\frak p^\nu$, where $\frak p$ is the prime ideal in $A$ defining $Y$ and
$\nu=Min\{n\in \Bbb N\mid e\leq $$n+r\choose r$$\}$ [Theorem 19].\medskip 
Note that
$Proj(G(R_\frak q) \otimes_{k(\frak p)}  K)$  is reduced if and only if 
there exists an open nonempty subset $U$ of $Y$ 
such that, for every closed point $x$ of $U$, the tangent cone to $X$ at $x$ 
is the union, as a set, of $e$ distinct linear spaces of $\Bbb A^{r+1}$
[Theorem 21]. Then
 (extending the definition given for curves)
 it is natural to say that in this case $Y$ 
is an {\sl ordinary multiple subvariety} of $X$.\medskip
Let $X=Spec R$ be a reduced non-normal $S_2$ variety. Assume
that the irreducible components of the non-normal locus of $X$ are ordinary multiple  
subvarieties $Y_i=Spec(R/\frak {q_i})$  of multiplicity  $e_i=e(R_\frak {q_i})>1$
($1\leq i\leq t$). Set
$emdim(R_\frak {q_i})=r_i+1$ and let $K_i$ be the algebraic closure of the residue field 
$k( \frak {q_i} )$. Suppose that the varieties $Y_i$ are nonsingular  and that
the points of $Proj(G(A_\frak {q_i})\otimes _{k( \frak {q_i} )} K_i)$ 
are in generic $e_i-1,e_i$  position, for any $i$.
Considering the following conditions:\medskip
(a) $X=Spec(R)$ is  Cohen-Macaulay, equimultiple of multiplicity
$e_i$ along $Y_i$ and has constant embedding dimension along $Y_i$, 
for any $i$;\medskip
(b) $X$ is normally flat along $Y_i$, for any $i$;\medskip
(c) $\frak b=\frak p_1^{\nu_1}\cap...\cap\frak p_t^{\nu_t}$;\medskip
then $(a) \Rightarrow (b) \Rightarrow (c)$ [Theorem 24]. \medskip
In particular
if, in (a), $X$ is a hypersurface, then the conductor
of $R$ is $\frak b=\frak p_1^{e_1-1}\cap...\cap\frak p_t^{e_t-1}$ [Corollary 26].

Throughout the paper all ring are supposed to be commutative, with identity and 
noetherian.\medskip
Let $A$ be a semilocal ring. If $\frak p$ is an  ideal of $A$, 
$G_\frak p (A)= \bigoplus_{n\geq 0} (\frak p^n/\frak p^{n+1})$ 
is the associated 
graded ring with respect to $\frak p$. By $G(A)$ 
we denote the associated graded ring 
with respect to the Jacobson radical $\frak J$ of $A$.
If $x\in A$, $x\neq 0$, $x\in \frak J^n- \frak J^{n+1}$, $n\in \Bbb N$, we say that $x$
has degree $n$ and the image $x^*\in \frak J^n/ \frak J^{n+1}$ of $x$ in $G(A)$
is said to be the initial form of $x$. If  $\frak a$ is an ideal of $A$, 
by $G(\frak a)$
we denote the ideal of $G(A)$ generated by all the inital forms of the elements
of $\frak a$.\medskip
 If $A$ is local with maximal ideal $\frak m$, 
$H^0(A,n)=dim_k(\frak m^n/ \frak m^{n+1})$, $n\in \Bbb N$, denotes the
 Hilbert function of $A$ and $e(A)$ the multiplicity of $A$ at $\frak m$. 
The embedding dimension $emdim(A)$ of $A$ is given by $H^0(A,1)$.
If $i\in \Bbb N$, the functions $H^i(A,n)$ are given by the relations 
$H^i(A,n)=$$\sum\nolimits_{j=0}^n$$ H^{i-1}(A,j)$

If $S=\bigoplus_{n\geq 0} S_n$ is a standard graded finitely generated
algebra over a field $k$,
of maximal homogeneous ideal $\frak n$,  
$H^0(S,n)=dim_k S_n=H^0(S_\frak n,n)$
 denotes the Hilbert function of $S$ and  $emdim(S)=H^0(S,1)=emdim(S_\frak n)$
the embedding dimension of $S$.  The multiplicity of
$S$ is  $e(S)=e(S_\frak n)$ .
  One has  $e(A)=e(G(A))$ 
and $emdim(A)=emdim(G(A))$.

If $B$ is any ring $dim(B)$ denotes the dimension
 of $B$.

\head 2. The one dimensional case.\endhead 
Let $B$ be a  local ring or a standard 
finitely generated graded $k$-algebra over a field $k$. 
Suppose $B$ has dimension one and
is Cohen-Macaulay. Set
$emdim(B)=r+1$ and $e(B)=e$.
It is well known that, for any $n\in\Bbb N$,
 $H^0(B,n) \leq Min\{e,$$n+r\choose r$$\}$ and if
$H^0(B,m)=e$ then $H^0(B,n)=e$ , for any $n\geq m$ [14].

\definition {1. Definition}The ring $B$ has {\sl maximal Hilbert function} if,
for every $n\in\Bbb N$,
$$H^0(B,n)=Min\{e,{n+r\choose r}\}$$ \enddefinition 
\definition {2. Definition} A set of points 
$X=\{P_1,...,P_e \}\subset\Bbb P^r$
is {\sl in generic position} (or the points $P_1,...,P_e$  are
  {\sl in generic position})
 [8, Definition 3.1] if the Hilbert function of
 its homogeneous coordinate ring
$R$ is maximal. The set $X$ is  {\sl in  generic $t$-position}, $t\leq e$,
 if
 every  $t$-subset of $X$ is in 
generic position (then generic $e$-position is generic position).\enddefinition

\remark {Remark} It is proved in [2, Theorem 4] that, for any $e$ and $r$, 
 \lq\lq generic position\rq\rq\ is an open nonempty condition.\endremark

\example {3. Example} It is easily seen that
any set of  
points of $\Bbb P^1$ is in generic t-position, for any $t$.\endexample
\example {4. Example} A set of  $n+r\choose r$ points in $\Bbb P^r$ $(n>0, r>0)$
is in generic position if and only if they do not lie on a hypersurface of degree
$n$ [8, Corollary 3.4], in particular six points in $\Bbb P^2$ are in generic position
if and only if they do not lie on a conic.\endexample

\proclaim{5. Theorem} Let B a one dimensional reduced local ring of maximal ideal $\frak m$. 
Let $K$ be the algebraic closure of the residue field $k( \frak m )=B / \frak m $.
Then:\medskip
(a) the Hilbert functions of $B$, $G(B)$ and  $G(B)\otimes_{k(\frak m)} K$
are the same, $e(B)=e(G(B))=e(G(B)\otimes_{k(\frak m)} K)$ and
 $emdim(B)=emdim(G(B))=emdim(G(B)\otimes_{k(\frak m)}  K)=r+1$;\medskip
(b) if  $Proj(G(B)\otimes _{k( \frak m )} K)$ is reduced and consists of points 
in generic position in $\Bbb P^r$ then $B$ has maximal Hilbert function and the
rings  $G(B)\otimes _{k( \frak m )}  K$ and $G(B)$ are reduced. 
\endproclaim

\demo{Proof} {\sl (a)} If $G(B)\otimes K$ is the $K$-vector space obtained extending to $K$
the field of scalars of the $k(\frak m)$-vector space $G(B)$, then the Hilbert functions of
$G(B)$ (that is of $B$) and of $G(B)\otimes K$ are the same. This implies the
equalities of the multiplicities and of the embedding dimensions.\medskip 
{\sl (b)}  If the points of $Proj(G(B)\otimes K)$ are in generic
position the Hilbert function of $G(B)\otimes K$ is maximal [Definitions 1 and 2].
Then by {\sl (a)}, the Hilbert function of  $B$ is maximal.
 This implies that $G(B)$ is Cohen-Macaulay [11, Theorem 3.2], that is,
there exists  $y^*\in \frak m/\frak m^2$ which is a non zero divisor of $G(B)$ (this
is easy to prove if $k(\frak m)$ is infinite and if $k(\frak m)$ is finite one can 
pass to the ring $R[U]_{\frak m R[U]}$, $U$ indeterminate).
Then, since a field extension is flat, $y^*\otimes 1$ is a non zero divisor of $G(B)\otimes  K$ which then is 
 Cohen-Macaulay and reduced, since $Proj(G(B)\otimes  K)$ is reduced.
Finally, by the flatness of $G(B)$ over $k(\frak m)$, we have 
$G(B)\subset G(B)\otimes K$ and $G(B)$ is reduced.
\qed\enddemo
\remark{Remark} If the ring $B$ of Theorem 5 is the local ring $A$ 
at a multiple point of a curve over an algebrically closed field $k$,
 then $K=k(\frak m)=k$ and 
$G(B)\otimes _{k( \frak m )} K=G(A)$. Hence in this case Theorem 5,(b) proves that,
if $Proj(G(A))$ is reduced and consists of points in generic position, then $G(A)$
is reduced which is Theorem 3.3 of [9]. In general the condition $Proj(G(A))$ reduced 
doesn't imply that $G(A)$ is reduced as has been shown with various examples 
in [9, Section 3, Example 1] and 
[7, Section 4]. Another example
(pointed out by A. De Paris ), which answers also to a question posed
in [3, Example 13], is the following.\endremark

\example{6. Example} Let $B=\Bbb C[g,tg,fg]$, $f=t^5-1$, $g=tf$ and $A$ be 
the local ring
of the curve $Spec(B)$ at $\frak m=(g,tg,fg)$. We have $e(A)=6$.
 Let $a_i$, $i=1,...,5$, be 
the fifth roots of the unity. It is proved in [7,Section 4] that 
 $Proj(G(A))$ consists of the points 
$P_i=(1,a_i,0) ,i=1,...,5, P_6=(1,0,-1)$ which lie on the conic $yz=0$ and then they
are not in generic position [Example 4].
But $G(A)$ is not reduced as we are going to show. Let 
$\frak n$ be the maximal homogeneous ideal of 
$G(A)_{red}=G(A)/nil(G(A))=\Bbb C[X,Y,Z]/\frak a=\Bbb C[x,y,z]$, 
where 
$\frak a=(Z,Y^5-X^5) \cap (Y,X+Z)$. 
There is
a natural surjective homomorphism 
 $\phi:\frak m^2/\frak m^3 \rightarrow\frak  n^2/\frak n^3$ 
given by $\phi(H(g,tg,fg))=H(x,y,z)$, $H(X,Y,Z)\in k[X,Y,Z]$
homogeneous of degree 2. Let $F(X,Y,Z)=Z^2+XZ$.Since $z^2+xz=0$,
if we show that $F(g,tg,fg)=fg(fg+g)\notin \frak m^3$
we have that $\phi$ is not injective and $G(A)$ is not reduced.
Now if $fg(fg+g)\in \frak m^3$  one  has 
$$fg(fg+g)=P(1,t,f)g^3+g^4h \tag *$$
where $h\in \Bbb C[t]$ and $P(X,Y,Z)$ is a homogeneous polynomial 
of degree $3$. But $fg(fg+g)=f(f+1)g^2=t^4g^3$. 
Substituting in (*) and dividing by $g^3$ we have $t^4-P(1,t,f)=gh$. But $t^4-P(1,t,f)=q+fr$
where $q,r\in \Bbb C[t]$ and $q$ is monic of degree $4$. Thus $q=gh-fr=f(th-r)$ which is impossible
because the degree of $f$ is $5$.\endexample

 We recall that the {\sl conductor} $\frak b$  of a ring  
$B$  in its normalization $\overline B$ is the ideal (of $B$ and $\overline B$)
$\frak b= \{ b\in$\ B$ \mid  \overline B b \subset B  \}$.

\proclaim{7. Theorem} Let $B$ be a one dimensional reduced local ring
 with finite 
normalization $\overline B$.
 Let $k(\frak m)=B/ \frak m$ be the residue field of $B$  
and $K$ be the algebraic closure of $k(\frak m)$. If $G(B) \otimes_{k(\frak m)} K$ is reduced 
then:
\medskip
(a) $G(B)$ is reduced, there is a natural 
immersion
$G(B) \subset G(\overline B)$ and $G(\overline B)$
is the normalization of $G(B)$;
\medskip
(b)
 there is a natural immersion 
$G(B) \otimes_{k(\frak m)} K \subset G(\overline B)\otimes_{k(\frak m)} K$ and 
$G(\overline B)\otimes_{k(\frak m)} K$
is the normalization of 
$G(B)\otimes_{k(\frak m)} K$;
\medskip
(c) If  $\nu=Min\{n\in\Bbb N\mid e=H(B,n)\}$,
 the ideal $G(\frak m)^\nu$ is contained in the conductor  of $G(B)$
in $G(\overline B)$, and if $B$ (that is $G(B)$) has maximal Hilbert function, then
$\nu=Min\{n\in\Bbb N\mid e\leq $$n+r\choose r$$\}$, 
where $e=e(B)$ and $r+1=emdim(B)$.
 
\endproclaim
\demo{Proof}{\sl (a)} We have $G(B)\subset G(B)\otimes K$ so if $G(B)\otimes K$
is reduced then so is $G(B)$. Then there is a natural 
immersion $G(B) \subset G(\overline B)$ [8, Proposition 2.2] and 
$G(\overline B)$ is the normalization of $G(B)$ [10, Proposition 1.7].\medskip
{\sl (b)}
By (a) and  by the flatness of $K$ over $k(\frak m)$ we have the immersion
$G(B)\otimes K \subset G(\overline B)\otimes K$. Furthermore  $G(\overline B)$
is the normalization of $G(B)$. From this we want to deduce that 
$G(\overline B)\otimes K$ 
is the normalization of 
$G(B)\otimes K$.
It is easily seen that   
$G(\overline B)\otimes K$ is integral over $G(B)\otimes K$ and 
that $G(\overline B)\otimes K$ is contained in the total ring of quotients of 
$G(B)\otimes K$. Since by assumption $G(B)\otimes K$ is reduced, 
 $G(\overline B)\otimes K$ is reduced. But, if $D= \overline B/ \frak m\overline B$,
$G(\overline B)\cong D[T]$ and $G(\overline B)\otimes K \cong (D\otimes K)[T]$.
Thus $D \otimes K$ is an artinian reduced ring, that is a direct sum of fields.
Hence $G(\overline B)\otimes K$ is normal.\medskip

{\sl (c)} If $\frak b$
is the conductor of $B$ in $\overline B$ it is well known 
([7, Lemma 2.12] and [12, Theorem 1.3]) that $\frak m^\nu
\subset\frak b$. Then $G(\frak m)^\nu \subset G(\frak b)$. But,
 $G(\frak b)$ 
is contained in the conductor of $G(B)$ in $G(\overline B)$.
In fact if $b^* \in G(\overline B)$ and $x^* \in G(\frak b)$ have degree respectively
$s$ and $t$ and are initial forms of
elements $b\in\frak J^s-\frak J^{s+1}$ and $x\in\frak J^t\cap\frak b-\frak J^{t+1}$
(where $\frak J$ is the Jacobson radical of $\overline B$),
then $bx\in \frak J^{s+t}\cap A=\frak m^{s+t}$ [8, Proposition 2.2] i.e. $b^*x^*\in G(A)$

Finally the last statement is an easy consequence of the definition of
maximal Hilbert function [Definition 1].\qed\enddemo

\proclaim{8. Theorem} Let $S$ be a standard graded finitely generated $k$-algebra, 
over an 
algebraically
closed field $k$, with maximal homogeneous ideal $\frak n$. Assume $S$ one dimensional
and reduced. Let $e=e(S)$ and 
$r+1=emdim(S)$. Then
$Proj(S)$ consists of $e$  points of $\Bbb P^r$. 
If these points are in generic $e-1,e$ position,
then the conductor of 
$S$ in its normalization $\overline S $ is
$\frak n^\nu$,where $\nu=Min\{n\mid e\leq $$n+r\choose r$$\}$.\endproclaim

\demo{Proof} [8, Proposition 3.5 and Theorem 4.3].\enddemo

\proclaim{9. Theorem} Let $A$ be the local ring of a curve at a singular point with
reduced tangent cone. Let $\frak m$ be the maximal ideal of $A$.
 If, for some integer $n$, $ G( \frak m^n)=G(\frak m)^n$
is the conductor of $G(A)$ in $G(\overline A)$, then
$\frak m^n$ is the conductor of $A$ in $\overline A$.\endproclaim
\demo{Proof} Let $\frak J$ be the Jacobson radical of $G(\overline A)$. Since by
assumption $ G(\frak m^n)$ is the conductor we have
$$ G( \frak m^n)= G( \frak m^n)G(\overline A)= G( \frak m^n\overline A)= G( \frak J^n)$$
 and then the result follows from [8, Theorem 2.3].\qed\enddemo

The following is the main result of [8] (see Theorem 4.4).

\proclaim{10. Theorem} Let $A$ be the local ring of a curve at 
a singular point, with 
maximal ideal $\frak m$. Let $e=e(A)$ and $emdim A=r+1$.
Assume that $Proj(G(A))$ is reduced  and consists of  points in generic e-1,e position .
Then $\frak m^\nu$ is the conductor of $A$, where 
$\nu=Min\{n\mid e\leq $$n+r\choose r$$\}$.
\endproclaim

\demo{Proof} By Theorem 5, (b) (see also the following Remark ) $G(A)$ is reduced.
Then, by Theorem 8, the conductor of $G(A)$ in its normalization $G(\overline A)$
is $G(\frak m)^\nu=G(\frak m^\nu)$ (note that $G(\frak m)$ is the maximal homogeneous 
ideal of $G(A)$).
Hence, by Theorem 9,
$\frak m^\nu$ is the conductor of $A$ in $\overline A$.\qed\enddemo

Theorem 10 can be extended to any one dimensional local ring in the following way:

\proclaim{11. Theorem} Let $B$ be a one dimensional reduced local  ring 
with finite 
normalization $\overline B$ and maximal ideal $\frak m$. Let $k(\frak m)=B/ \frak m$ be 
the residue field of $B$  
and $K$ be the algebraic closure of $k(\frak m)$.
Set $e=e(B)$ and $emdim(B)=r+1$.
If $Proj(G(B)\otimes_{k(\frak m)} K)$ is reduced and consists of points in generic $e-1,e$
 position
then the conductor  of $B$ in $\overline B$ is $\frak m^\nu$ where 
 $\nu=Min\{n\mid e\leq $$n+r\choose r$$\}$.
\endproclaim

\demo{Proof} If $Proj(G(B)\otimes K)$ is reduced and consists of points in generic 
$e-1,e$ position then 
$G(B)\otimes K$ is reduced [Theorem 5, (b)], its
normalization is $G(\overline B)\otimes  K$ [Theorem 7,(b)], 
and its conductor is  $(G(\frak m)\otimes K)^\nu= 
 G(\frak m)^\nu \otimes K$ [Theorem 8 and Theorem 5, (a)]. 
 Furthermore, if $\frak b$ is the conductor of $G(B)$, 
then
$\frak b \otimes K$ is contained in the conductor 
of $(G(B)\otimes K)$.
In fact if
$y\in\frak b$ , $b\in G(\overline B)$ and $c,c'\in K$,
$(b\otimes c)(y \otimes c')=by\otimes cc'\in  G(B)\otimes K$.
But, by Theorem 7,(c), $G(\frak m)^\nu \subset \frak b $ and then, by flatness
$G(\frak m)^\nu \otimes K \subset \frak b \otimes K $.
Thus
$$G(\frak m)^\nu \otimes K = \frak b \otimes K$$
 But, again by flatness, $$dim_k(\frak b /G(\frak m)^\nu)\otimes K=
dim_k((\frak b \otimes K)/(G(\frak m)^\nu\otimes K))=0$$
that is $\frak b=G(\frak m)^\nu$ and, by Theorem 9,  
the conductor of $B$ is $\frak m^\nu$.\qed\enddemo

\remark{Remark}Theorem 10 is the main result of a more general statement
on the conductor of a curve at an ordinary singularity [8, Theorem 4.4]. This statement
can be easily extended to  one dimensional rings
in the same way of Theorem 11.\endremark

\head 3. The codimension one case.\endhead

We need some preliminaries.\medskip

\definition{12. Definition} Le $\frak p$ be a prime ideal of a local ring $A$.
$A$ is {\sl normally flat along $\frak p$ } if $\frak p^n/\frak p^{n+1}$
is flat over $A/\frak p$, for all $n\geq 0$. Note that $A$ is normally flat
along $\frak p$ if nad only if $G_\frak p(A)$ is free over $A/\frak p$.\enddefinition

\proclaim{13. Theorem} Let  $\frak p$ be a prime ideal of a local ring $A$ such that
  $A/\frak p$ is regular and $dim(A/\frak p)=d \geq 1$.
Then:\medskip
 $H^0(A,n) \geq H^d(A_\frak p,n)$, for any n and equality holds if and only if
$A$ is normally flat along $\frak p$.\medskip

In particular $$emdim(A) \geq emdim(A_\frak p)+dim(A/\frak p), e(A) \geq e(A_\frak p)$$
and if $A$ is normally flat along $\frak p$, 
then $$emdim(A)=emdim(A_\frak p)+dim(A/\frak p), 
e(A)=e(A_\frak p)$$.\endproclaim
\demo{Proof} See [5, Theorem (22.24) and Proposition (30.2)].\qed\enddemo

\proclaim{14. Theorem} Let $\frak p$ be a prime ideal of codimension one of 
a Cohen-Macaulay local ring $A$ such that $A/\frak p$ is regular.
 Assume  $e(A)=e(A_\frak p)=e$ and 
$emdim(A)=emdim(A_\frak p)+dim(A/\frak p)$.
If $H^0(A_\frak p,n)$ is maximal 
 then $A$ is normally flat along $\frak p$
and $G(A)$ is Cohen-Macaulay.\endproclaim
\demo{Proof} Set $dim(A/\frak p)=d$, $e(A_\frak p)=e$ and $emdim(A_\frak p)=r+1$.
Since $A$ is Cohen-Macaulay there is an $A$-sequence ${x_1,...,x_d}$ 
of elements  of degree 1 (this can always be arranged if the residue field of $A$ is infinite and
if it is finite, by passing, if necessary, to the ring
$A[U]_{\frak m A[U]}$, $U$ indeterminate) and such that the ring
 $B=A/(x_1,...,x_d)$ is Cohen-Macaulay. Then, by the assumption, we have $e(B)=e(A)=e(A_\frak p)=e$ , 
$emdim(B)=emdim(A)-d=emdim(A_\frak p)=r+1$. Furthermore $B$ and $A_\frak p$ are 
one dimensional and then   $H^0(B,n) \leq Min\{e,$$n+r\choose r$$\}=H^0(A_\frak p,n)$.
[Definition 1 and preceding comments]
Hence $H^d(B,n) \leq H^d(A_\frak p,n)$. But 
$H^d(B,n)\geq H^0(A,n)$ [11, Theorem 3.1 (note that in the 
statement there is a misprint, namely the exponent $s$ of the Hilbert function $H^s$ 
should be $s+1$)]. Furthermore $H^0(A,n) \geq H^d(A_\frak p,n)$ by Theorem 13.
We have then $H^d(B,n)=H^0(A,n)=H^d(A_\frak p,n)$. Hence $A$ is normally flat along $\frak p$
[Theorem 13] and  the initial forms $x_1^*,...,x_d^*$ form a regular sequence of
$G(A)$ [11, Theorem 3.1]. Furthermore, since $H^0(B,n)=H^0(A_\frak p,n)$, the 
one dimensional ring $B$ has maximal Hilbert function and then $G(B)$ is
 Cohen-Macaulay [11, Theorem 3.2]. But $G(B)\cong G(A)/(x_1^*,...,x_d^*)$
 [14, Chapter 2, Lemma 3.2] hence 
$G(A)$ is Cohen-Macaulay .\qed\enddemo

\proclaim{15. Theorem} Let $\frak p$ be a prime ideal of codimension one of 
 a reduced local ring $A$ with finite normalization $\overline A$. Let
$e(A_\frak p)=e>1$, $emdim(A_\frak p)=r+1$ and  denote by
$K$ the algebraic closure of the residue field $k(\frak p)=A_\frak p/\frak pA_\frak p$
of $A$ at $\frak p$.
Let        $\frak b$ be the conductor of $A$
in  $\overline A$.
Assume $A/\frak p$  
regular, $\sqrt{ \frak b}=\frak p$ and $Proj(G(A_\frak p) \otimes_{k(\frak p)}  K)$ 
 reduced
with   points in generic $e-1,e$ position.\medskip
Given the following   conditions:\medskip 

(a) $A$ is Cohen-Macaulay and $emdim(A)=emdim(A_\frak p)+dim(A/\frak p)$, $e(A)=e(A_\frak p)$;  \medskip

(b)  $A$ is $S_2$   and   normally flat 
along $\frak p$;\medskip

(c) $\frak b$ is primary and $A$ is normally flat 
along $\frak p$;\medskip

(d) $\frak b=\frak p^\nu$, where $\nu=Min\{n\in \Bbb N\mid e\leq $$n+r\choose r$$\}$;
\medskip
then $(a) \Rightarrow (b) \Rightarrow (c) \Rightarrow (d)$.\medskip

\endproclaim

\demo{Proof} $(a) \Rightarrow (b)$ A  Cohen-Macaulay ring is $S_2$ [6, 17.I]. By 
Theorem 5,(b) if $Proj(G(A_\frak p) \otimes_{k(\frak p)} K)$
is reduced and consists of points in generic position, the Hilbert function of $A_\frak p$
is maximal and, by Theorem 14, $A$ is normally flat along $\frak p$.\medskip 
 $(b) \Rightarrow (c)$ If $A$ is  $S_2$ then its conductor $\frak b$ is unmixed
 [4,Lemma 7.4] and then primary,since $\sqrt{ \frak b}=\frak p$ .\medskip
 $(c) \Rightarrow (d)$ $\frak b A _\frak p$ is the conductor of $A _\frak p$ in
its normalization $\overline A _\frak p$ and then, by Theorem 11,
 $\frak b A _\frak p$=$
(\frak p A_\frak p)^\nu$=$\frak {p^\nu}A_\frak p$.
 Furthermore, if $A$ is normally flat along $\frak p$ then $\frak p^\nu$
 is primary [13, Proposition 1.1] and then
$$\frak b=A \cap \frak b A _\frak p=A \cap\frak {p^\nu}A_\frak p=\frak p^\nu$$ 
\qed\enddemo

\proclaim{16. Corollary} Suppose $\frak p$ is a prime ideal of codimension one of a 
Cohen-Macaulay reduced
local  ring $A$ with finite normalization $\overline A$.
Denote by
$K$ the algebraic closure of the residue field 
$k(\frak p)=A_\frak p/\frak pA_\frak p$.Let $\frak b$ be the conductor of $A$
in  $\overline A$.
Assume $A/\frak p$  
regular, $\sqrt{ \frak b}=\frak p$ and $Proj(G(A_\frak p) \otimes_{k(\frak p)}  K)$ 
 reduced. If $emdim(A)=dim(A)+1$ and $e(A)=e(A_\frak p)=e$
then $\frak b=\frak p^{e-1}$.\endproclaim
\demo{Proof}Let $emdim(A)=dim(A)+1$. By assumption $dim(A)=dim(A/\frak p)+1$
and, by Theorem 13, $emdim(A) \geq emdim(A_\frak p)+dim(A/\frak p)$. Then
$emdim(A_\frak p)=2$ ($emdim(A_\frak p)>1$ since $e(A_\frak p)>1$) and
 $emdim(A)=emdim(A_\frak p)+dim(A/\frak p)$. Hence 
$emdim(G(A_\frak p) \otimes_{k(\frak p)} K)=2$ [Theorem 5,a)] that is
$Proj(G(A_\frak p) \otimes_{k(\frak p)}  K)\subset\Bbb P^1$ and then its points
are always in generic $e-1,e$ position [Example 3]. Hence the claim follows from
Theorem 15, $(a) \Rightarrow (d)$ (if $r=1$ one has $\nu=e-1$).\qed\enddemo

In the following we want to apply the previous results to the geometric case.  
\proclaim {Standing notation} From now on $X=Spec(R)$ is a reduced variety over an algebraically closed
field $k$ and $Y=Spec(R/\frak q)$ is an irreducible   codimension one subvariety $Y$ 
of $X$ (i.e. $\frak q$ is a prime ideal of codimension one in $R$). By $K$ we  
denote 
the algebraic closure of the residue field
 $k( \frak q )$ of $R$ at $\frak q$. We set $emdim(R_\frak q)=r+1$ and 
$e(R_\frak q)=e$.\endproclaim

\definition {17. Definition} Let $x$ be any closed point of $Y$ and $A$ 
be the local ring of $X$ 
at $x$. Let $\frak p=\frak qA$ be 
the prime ideal in
 $A$ defining the subvariety $Y$. Then:\medskip
(i) $Y$ is a {\sl multiple subvariety} of $X$ of
multiplicity $e$ if $e(R_\frak q)=e>1$;\medskip 
(ii) $X$ is {\sl normally flat along $Y$ at $x$} if 
$A$ is normally flat along $\frak p$;
\medskip
(iii) $X$ is {\sl normally flat along} $Y$ if it is so at every point of $Y$;\medskip
(iv) $Y$ is 
{\sl nonsingular
at} $x$ if $A/\frak p$ is regular.\enddefinition
\remark{Remark} With the notations of Definition 17 we have $A_\frak p=R_\frak q$ 
for any closed point $x$ of $Y$.\endremark

\proclaim{18. Theorem} Let  $Y$ be a 
multiple subvariety of $X$.
There exists an open nonempty subset $U$ of $Y$ such that,
 for every closed point $x$
of $U$,  $Y$ is nonsingular at $x$ and
 $X$ is normally flat along $Y$ at $x$.\endproclaim
\demo{Proof} It is well known  that 
the nonsingular points of $Y$
form an open nonempty set
and that  $X$ is normally flat along $Y$ at the points of an open nonempty
subset of $Y$ [5, Corollary (24.5)].\qed\enddemo

\proclaim{19. Theorem} Set $emdim(R_\frak q)=r+1$. Suppose that $Y$ is a 
multiple subvariety of $X$ of multiplicity $e$. If  $Proj(G(R_\frak q)\otimes _{k( \frak q )} K)$ 
is reduced and consists 
of points in
generic $e-1,e$ position then there exists an open nonempty subset $U$ of $Y$ 
such that,
for every closed point $x$ of $U$,  
the conductor $\frak b$ of the local ring $A$ of $X$ at $x$
 is primary and equal to $\frak p^\nu$ where $\frak p$ is the prime ideal in $A$ 
defining $Y$ and 
$\nu=Min\{n\mid e\leq $$n+r\choose r$$\}$.\endproclaim 
\demo{Proof} By the Remark to Definition 17, if $x$ is any point of $X$ and $A$ is
the corresponding local ring,
 one has $e(A_\frak p)=e(R_\frak q)$,  
hence, by assumption, $e(A_\frak p)>1$ that is, since $dim(A_\frak p)=1$,
$A_\frak p$ is not normal.
Then it is well known that $\frak b$ is contained in $\frak p$ . 
But the conductor contains a nonzero divisor.
Thus the codimension one ideal $\frak p$ is 
a minimal prime  of $\frak b$ . Now
it is easily shown that there exists an open nonempty subset $U_1$ of $Y$ such that,
if $x$ is a closed point of $U_1$, 
$\frak p$ is the unique prime ideal
of  $A$  associated to 
$\frak b$. Then $\frak b$ is primary and $\sqrt{ \frak b}=\frak p$.
 Furthermore, by Theorem 18, there exists an open nonempty subset $U_2$ of $Y$
such that
 for every closed point $x$
of $U_2$, $Y$ is nonsingular at $x$ and
 $X$ is normally flat along $Y$ at $x$. Then, if we apply 
Theorem 15, $(c) \Rightarrow (d)$ to the local rings of the points of $U=U_1\cap U_2$, 
we have the claim.
\qed\enddemo 
We want now to characterize geometrically the condition:  
$Proj(G(R_\frak q)\otimes _{k( \frak q)} K)$ is reduced. We need the following
 preliminary
general result.\medskip

\proclaim{20. Theorem} Let $A$ be the local ring of $X$ at a closed point $x$ 
of $Y$
and $\frak p$ be the prime ideal defining $Y$ in $A$.
 Set $dim(A/\frak p)=d$. Suppose $Y$ is nonsingular at  $x$ and $X$ is normally flat along $Y$
at $x$. Then there is an
isomorphism of graded k-algebras
  
$$G(A) \cong (G_\frak p (A) \otimes _{A/ \frak p} k)[T_1,...,T_d]$$ \endproclaim
\demo{Proof} See [5,Corollary (21.11)].\qed\enddemo

Assume that $Y$ is a multiple subvariety of $X$ and $A$ is the local ring of $X$ 
at a closed point $x$ of $Y$.
$A/\frak p$ is the local ring of  $Y$  at $x$.
Then $Spec(G(A))$ is the  tangent cone to $X$ at $x$ and
$Spec(G(A/\frak p))$ is the tangent cone to $Y$ at $x$. Furthermore, there is a natural 
surjective homomorphism $G(A) \rightarrow G(A/\frak p)$ and then $Spec(G(A/\frak p))$
naturally embeds in $Spec(G(A))$. In this setting we have the following result:

\proclaim{21. Theorem} Let $K$ be 
the algebraic closure of the residue field
 $k( \frak q)$ of $R$ in $\frak q$. Then 
 $Proj(G(R_\frak q)\otimes _{k( \frak q)} K)$ is reduced (that
is consists 
of $e$ points) if and only if there exists an open nonempty subset $U$ of $Y$ 
such that, for every closed point $x$ of $U$, the tangent cone to $X$ at $x$ 
is the union, as a set, of $e$ distinct linear varieties,whose intersection
is the tangent cone to $Y$ at $x$.
\endproclaim 

\demo{Proof} If $Proj(G(R_\frak q)\otimes _{k( \frak q )} K)$ has $e$ points  by [1, Theorem 1.1]
there exists an open nonempty subset $U_1$ of $Y$ such that, for every closed point $x$
of $U_1$, the
tangent cone  $Spec(G(A))$ to $X$ at $x$ , has $e$ irreducible components,
 that is 
$G(A)$
has $e$ minimal primes. Furthermore, by Theorems 18 and 20, there exists an open nonempty
subset  $U_2$ of $Y$ such that, for every closed point $x$ of $U_2$,  
 $G(A)$ is a polynomial ring over 
$G_\frak p (A) \otimes _{A/ \frak p} k$. Then if $x$ is a point of $U_2$
 the minimal primes of $G(A)$ are extensions of
 the minimal primes of the one dimensional finitely generated graded $k$-algebra
 $G_\frak p (A) \otimes _{A/ \frak p} k$, hence they are generated by linear forms. 
Since, at every closed point $x\in U=U_1 \cap U_2$ 
 the tangent cone to $Y$  is contained in the tangent cone to 
$X$,          we have the claim.
Vice versa if, at every point of an open nonempty subset of $Y$, the tangent cone 
to $X$ has $e$ irreducible
components, again by [1, Theorem 1.1], also the zero dimensional scheme 
$Proj(G(A_\frak p)\otimes _{k( \frak p )} K)$
has $e$ irreducible components, 
that is $e$
points, and then, since by Theorem 5,(a) $e(G(A_\frak p)\otimes _{k( \frak p )} K)=
e(G(A_\frak p))=e(A_\frak p)=e$, is reduced.\qed\enddemo

\remark{Remark} Theorem 21 extends Corollary 1.4 of [1] which states
the only if part  in the particular case in which 
$G(R_\frak q)\otimes _{k( \frak q)} K$ is reduced.
 \endremark\medskip
\definition{22. Definition} Let  $Y$ be a multiple subvariety of $X$. If 
$Proj(G(R_\frak q)\otimes _{k(\frak q)} K)$ is reduced 
$Y$ is said to be an {\sl ordinary multiple subvariety} of $X$.\enddefinition 

\remark{Remark} If $Y$ is a multiple point of a curve $X$ the previous definition 
agrees the one given in  [9, Lemma-Definition 2.1].\endremark

\definition{23. Definition} If $A$
is the local ring of the variety $X=Spec(R)$ at a closed point $x$, then $e(A)$ and $emdim(A)$
are called respectively the {\sl multiplicity} and the {\sl embedding dimension}
{\sl of $X$ at $x$}. Let $Y=Spec(R/\frak q)$ be a multiple subvariety of $X$. Then\medskip  
(i) $X$ is {\sl equimultiple along} $Y$ if $X$ has the same multiplicity $e(R_\frak q)$
at every point of $Y$; \medskip
(ii) $X$ has {\sl constant embedding dimension along} $Y$ 
if $X$ has the same embedding
dimension $emdim(R_\frak q)+dim(R/\frak q)$ at every point of $Y$.\enddefinition

 \remark {Remark} By Theorem 18 there exists one point  $x$ of $Y$ for which
 $Y$ is nonsingular at $x$ and  $X$ is normally flat along $Y$ at $x$,
then by Theorem 13 the multiplicity of $X$ at $x$ is $e(R_\frak q)$ and the 
embedding dimension of $X$ at $x$ is $emdim(R_\frak q)+dim(R/\frak q)$ (see also
the Remark to Definition 17); hence if
$X$ is equimultiple and has constant embedding dimension along $Y$ these 
have to be the multiplicity and embedding dimension at all points of $Y$.\endremark

\medskip
It is well known that, if $X=Spec(R)$ is a reduced non-normal $S_2$  variety,
the non-normal locus of $X$ (that is the set of closed points in $X$ whose local ring
is not normal) consists of the union of irreducible  subvarieties 
of codimension one, $Y_i=Spec(R/\frak {q_i})$, ($1\leq i\leq t$). Furthermore the
$\frak q_i$ are the primes associated to  the conductor of $R$,
that is the primes for 
which the one dimensional rings $R_\frak {q_i}$ are not normal
 or equivalently not regular.
\proclaim{24. Theorem} Let $X=Spec R$ be a reduced non-normal $S_2$ variety. Assume
that the irreducible components of the non-normal locus of $X$ are  ordinary multiple  
subvarieties $Y_i=Spec(R/\frak {q_i})$
($1\leq i\leq t$)   of multiplicity  $e_i=e(R_\frak {q_i})>1$. Set
$emdim(R_\frak {q_i})=r_i+1$ and let $K_i$ be the algebraic closure of the residue field 
$k( \frak {q_i} )$. If $Proj(G(R_\frak {q_i})\otimes _{k( \frak {q_i})} K_i)$ 
has points in generic $e_i-1,e_i$  position, for any $i$, 
then the conductor of $R$ is 
$\frak b=\frak q_1^{(\nu_1)}\cap...\cap\frak q_t^{(\nu_t)}$
where $\nu_i=Min\{n\mid e_i\leq $$n+r_i\choose r_i$$\}$ 
and $\frak q_i^{(\nu_i)}$ denotes the $\nu_i-th$ symbolic power of $\frak q_i$. 
Furthermore assume that the  $Y_i$ are  nonsingular varieties, for any $i$, 
and consider the following conditions:\medskip
(a) $X=Spec(R)$ is  Cohen-Macaulay, equimultiple 
 along $Y_i$ and has constant embedding dimension along $Y_i$ , 
for any $i$;\medskip
(b) $X$ is normally flat along $Y_i$, for any $i$;\medskip
(c) $\frak b=\frak q_1^{\nu_1}\cap...\cap\frak q_t^{\nu_t}$;\medskip
then $(a) \Rightarrow (b) \Rightarrow (c)$. \endproclaim
\demo{Proof} Let $\frak b=\frak a_1\cap...\cap\frak a_t$ be a minimal primary
decomposition, $\sqrt{ \frak a_i}=\frak q_i$ ($1\leq i\leq t$). By Theorem 11,
 for any $i=1,...,t$, 
$\frak {a_i} R_\frak {q_i}=\frak b R_\frak {q_i}=\frak {q_i}^{\nu_i} R_\frak {q_i}= 
\frak {q_i}^{(\nu_i)} R_\frak {q_i}$.
Then 
$$\frak {a_i}=R \cap \frak {a_i} R _\frak {q_i}=R \cap\frak {q_i^{\nu_i}}R_\frak {q_i}=
\frak q_i^{(\nu_i)}$$
which is the first claim. 
The implications of the second claim follow easily from the analogous local
implications of Theorem 15.
\qed\enddemo
If we apply Theorem 24 to the particular case $t=1$ we get that, 
under the assumption $(a)$,
Theorem 19 holds for any point $x$ of the subvariety $Y$ (not only on an open
subset of $Y$). 
\proclaim{25. Corollary } Let $X=Spec R$ be a reduced non-normal Cohen-Macaulay  variety.
Assume that the non-normal locus of $X$ is an ordinary multiple nonsingular 
subvariety $Y=Spec(R/\frak q)$ of multiplicity $e=e(R_\frak q)>1$ and that $X$
 is equimultiple and has constant embedding
dimension along $Y$. Then, for every  point $x$ of $Y$,  
the conductor $\frak b$ of the local ring $A$ of $X$ at $x$
 is primary and equal to $\frak p^\nu$ where $\frak p$ is the prime ideal in $A$ 
defining $Y$ and 
$\nu=Min\{n\mid e\leq $$n+r\choose r$$\}$\endproclaim 
\demo{Proof} Setting $t=1$ in Theorem 24 we have that 
$\frak b R_\frak q=\frak {q^\nu} R_\frak q=(\frak q R_\frak q)^\nu=\frak p^{\nu}$
is the conductor of $A$.\qed\enddemo

\proclaim{26. Corollary} Let $H$ be a non-normal hypersurface. Assume that the 
irreducible
components of the non-normal locus of $X$ are  ordinary multiple  
subvarieties $Y_i=Spec(R/\frak {q_i})$   
 of multiplicity  $e_i=e(R_\frak {q_i})>1$ ($1\leq i\leq t$).
 If the $Y_i$ are nonsingular varieties 
and $X$ is equimultiple 
along $Y_i$,  
for any $i$, then the conductor
of $R$ is $\frak b=\frak q_1^{e_1-1}\cap...\cap\frak q_t^{e_t-1}$.\endproclaim

\demo{Proof}If $H=Spec(R)$ is a hypersurface, then $R$ is Cohen-Macaulay and
its embedding dimension at any point is constant. Furthermore  
$Proj(G(R_\frak {q_i}) \otimes_{k(\frak {q_i})}  K_i)\subset\Bbb P^1$ has points
 in generic $e_i-1,e_i$
position [Example 3]. Then the claim follows by Theorem 24, $(a)\Rightarrow (c)$,
in which $r_i=1$.\qed\enddemo

\example{27. Example} Let $H=Spec R$, $R=\Bbb C[T_0,...,T_r]/(L_1...L_n)$ 
($r\geq 2, n\geq 2$)
 be the union of 
$n$ hyperplanes of $\Bbb P^r$. Let $\frak q_i$ ($1\leq i\leq t$) denote the codimension
one primes of $R$
of the form $(\overline{L}_p,\overline{L}_q)$, $p,q\in \{1,...,n\}$, $p\not=q$, and 
  $Y_1,...,Y_t$ be the corresponding linear varieties  (that is the irreducible
components of the non-normal locus of $X$). Then 
 $e(R_{\frak q_i})=e_i$ is the number of hyperplanes passing through $Y_i$. 
Furthermore since $Y_i$ are linear it is easily shown that $\frak q_i^m=\frak q_i^{(m)}$, 
for any $m$. Hence from Theorems 21 and 24, it follows that 
 the conductor of $R$ is $\frak b=\frak q_1^{e_1-1}\cap...\cap\frak q_t^{e_t-1}$.\endexample  

\example{28. Example} Let $R=\Bbb C[X,Y,Z]/(XY^n-Z^n)=\Bbb C[x,y,z]$ ($n\geq 2$) .
The non-normal locus of the hypersurface $H=SpecR$ is the line 
$L:y=0,z=0$ of multiplicity $n$ , that is $e(R_\frak q)=n$, where $\frak q=(y,z)$.
Moreover it is easily seen that $H$ is equimultiple along $L$. Furthermore, if
$a\in\Bbb C$, $a\neq 0$, the tangent cone,
 at the  point $(a,0,0)$ of $L$,  
consists, as a set, of the $n$ distinct
planes $z=by$ , where $b^n=a$ , 
and then, by Theorem 21 and Corollary 26 and, the  conductor of $R$ in 
$\overline R$ is $(y,z)^{n-1}$.\endexample
\remark {Acknowledgements} The author would like to thank the referee for his
comments and suggestions.\endremark

\enddocument